\documentclass[10pt,reqno]{amsart}
\usepackage{amscd,amsmath,amssymb,amsthm}
\title{Nilpotent extensions of minimal homeomorphisms}
\author{Gernot Greschonig}
\author{Ulrich Hab\"ock}
\address{Faculty of Mathematics, University of Vienna, Nordbergstra\ss e 15, A-1090 Vienna, Austria} \email{\{gernot.greschonig\}\{ulrich.haboeck\}@univie.ac.at}
\date{}
\newtheorem{theo}{Theorem}[section]
\newtheorem{cor}[theo]{Corollary}
\newtheorem{lem}[theo]{Lemma}
\newtheorem{prop}[theo]{Proposition}
\theoremstyle{definition}
\newtheorem{defi}[theo]{Definition}
\theoremstyle{remark}
\newtheorem{rem}[theo]{Remark}
\newtheorem{rems}[theo]{Remarks}
\newtheorem{rem*}[]{Remark}

\newcommand{\comment}[1]{}

\newcommand{\Z}{\mathbb Z}

\newcommand{\R}{\mathbb R}

\newcommand{\mc}{\mathcal}

\DeclareMathOperator{\Stab}{\it{Stab}}
\thanks{The authors were partially supported by the FWF research project P16004-MAT}

\begin{document}\allowdisplaybreaks\frenchspacing

\begin{abstract}
In this paper we study topological cocycles for minimal homeomorphisms on a compact metric space.
We introduce a notion of an essential range for topological cocycles with values in a locally compact group, and we show that this notion coincides with the well known topological essential range if the group is abelian.
We define then a regularity condition for cocycles and prove several results on the essential ranges and the orbit closures of the skew product of regular cocycles.
Furthermore we show that recurrent cocycles for a minimal rotation on a locally connected compact group are always regular, supposed that their ranges are in a nilpotent group, and then their essential ranges are almost connected.
\end{abstract}

\maketitle

\section{Introduction}
\label{s:intro}
The aim of this paper is to study the dynamical behaviour of topological cocycles with values in a nonabelian locally compact second countable (l.c.s.) group  $G$, and throughout the paper these cocycles will be defined over a minimal homeomorphism on a compact metric space.
In the paper \cite{A} Atkinson introduced the notion of essential values for continuous $\R^n$-valued cocycles and studied their skew product extensions.
Later his results were extended to abelian groups by Lema\'nczyk and Mentzen (see \cite{LM} and also \cite{M1}).
But if one uses this notion of an essential range naively in the case of a non-abelian group (as has been done in \cite{LM}), it turns out that one cannot obtain much insight into the dynamical system.
In fact, in \cite{M2} examples of two cohomologous cocycles are given, one of which has a non-trivial essential range while the other one has a trivial essential range.

Let $G$ be a l.c.s. group with identity $\mathbf{1}_G$, let $T$ be a minimal homeomorphism of a compact metric space $(X,\delta)$ and let $f\colon X\longrightarrow G$ be a continuous map.
We define then a map $f\colon \mathbb{Z}\times X \longrightarrow G$ by
\begin{equation}
  \label{eq:cocycle}
  f(n,x)=
        \begin{cases}
          f(T^{n-1} x)\cdots f(T x)\cdot f(x) & \textup{if}\enspace n\geq 1
        \\
        \mathbf{1}_G & \textup{if}\enspace n=0
        \\
        f(-n,T^n x)^{-1} & \textup{if}\enspace n < 0.
    \end{cases}
\end{equation}
This map satisfies that
\begin{equation}\label{eq:cocycle2}
  f(k,T^l x)\cdot f(l,x)=f(k+l,x)
\end{equation}
for all integers $k,l$ and every $x\in X$, and is thus a \emph{cocycle} of the $\mathbb{Z}$-action on $X$ defined by $(m,x)\longmapsto T^m x$.

The \emph{skew product} transformation of $f$ is the homeomorphism on $X\times G$ defined by
\begin{equation}
  \mathbf{T}_f(x,g)=\big(T x, f(x)\cdot g\big),
\end{equation}
and this transformation related to the map $f(n,x)$ by the equality that
\begin{equation*}
  \mathbf{T}_f^n(x,g)=\big(T^n x, f(n,x)\cdot g\big)\enspace\textup{for all}\enspace n\in\mathbb{Z}.
\end{equation*}
Throughout this paper we denote the orbit closure of $(x,g)\in X\times G$ under $\mathbf{T}_f$ by
\begin{equation*}
\bar{\textup{O}}_f(x,g)=\overline{\big\{\mathbf{T}_f^n(x,g):n\in\Z\big\}}.
\end{equation*}

A cocycle is called \emph{topologically recurrent} if for every open neighbourhood $U$ of $\mathbf{1}_G$ and every open set $\mc O\subseteq X$ there is an integer $n\neq 0$ so that
\begin{equation*}
T^{-n}\mc O \cap \mc O \cap \big\{x: f(n,x)\in U\big\}\neq\varnothing.
\end{equation*}
This property is equivalent to the \emph{topological conservativity} (regional recurrence in the terminology of \cite{G-H}) of the skew product $\mathbf{T}_f$, i.e. for every open set $\mc O'\subseteq X\times G$ there is an integer $n\neq 0$ so that $\mathbf{T}_f^n(\mc O')\cap\mc O'\neq\varnothing$.

For nonabelian extensions the following ``local'' notion of an essential range is suitable to study the dynamics of the system:

\begin{defi}
A group element $g\in G$ belongs to the \emph{essential range} of the cocycle $f$ at $x\in X$, which is denoted by $E_x(f)$, if for every open neighbourhood $U$ of $g$ and every open neighbourhood $\mc O$ of $x$ there exists an integer $n\neq 0$ so that
\begin{equation}\label{eq:ev}
T^{-n}\mc O \cap \mc O \cap \big\{x: f(n,x)\in U\big\}\neq \varnothing.
\end{equation}
It is obvious from the definition that the essential range is a closed subset of $G$, and from the cocycle equality (\ref{eq:cocycle2}) it follows that $f(-n,T^n x)=f(n,x)^{-1}$ and hence $E_x(f)$ is symmetric, i.e. $E_x(f)=E_x(f)^{-1}$.
If $f$ is topologically recurrent then the identity is in $E_x(f)$ for every $x\in X$.
\end{defi}

\begin{rems}
(i) In difference to the definition above the topological essential range introduced in \cite{LM}, denoted by $E(f)$, does not refer to some $x\in X$.
For any $g\in E(f)$ the equation (\ref{eq:ev}) must be fulfilled for every open neighbourhood $V$ of $g$ and \emph{every open} set $\mc O\subseteq X$, and thus
\begin{equation*}
E(f)=\bigcap_{x\in X} E_x(f).
\end{equation*}

(ii) Assume that $x_n\to x$ and $g_n\to g$ are two convergent sequences in $X$ and $G$ respectively, and assume that $g_n\in E_{x_n}(f)$ for all positive integers $n$.
Then the definition of the essential range shows that $g\in E_x(f)$.
\end{rems}

An important property of the local essential ranges is that they are always conjugate for points belonging to the same $T$-orbit in $X$:
\begin{lem}\label{lm:E_orbit}
For every $x\in X$ and every integer $n$ we have the following equality:
\begin{equation}\label{eq:ER_ORBIT}
E_{T^n x}(f)=f(n,x)\cdot E_x(f)\cdot f(n,x)^{-1}
\end{equation}
\end{lem}

\begin{proof}
Let $U$ be an open neighbourhood of $g\in f(n,x)\cdot E_x(f)\cdot f(n,x)^{-1}$ and $\mc O$ an open neighbourhood of $T^n x$.
We set $h=f(n,x)^{-1}\cdot g \cdot f(n,x)\in E_x(f)$ and choose a symmetric open neighbourhood $V$ of $\mathbf{1}_G$ so that $f(n,x)\cdot V h V^2 \cdot f(n,x)^{-1}\subseteq U$.
We choose then an open neighbourhood $\mc O'$ of $x$ so that $T^n \mc O'\subseteq \mc O$ and $f(n,\mc O')\subseteq f(n,x)\cdot V$, and as $h\in E_x(f)$ we can find an integer $k\neq 0$ and $y\in X$ with $y\in T^{-k} \mc O' \cap \mc O' \cap \big\{x: f(k,x)\in hV \big\}$.
It follows that $T^n y\in T^{-k}\mc O \cap \mc O$ while
\begin{equation*}
f(k,T^n y)=f(n,T^k y)\cdot f(k,y)\cdot f(n,y)^{-1}\in f(n,x)\cdot V h V^2\cdot f(n,x)^{-1}\subseteq U,
\end{equation*}
and as $U$ and $\mc O$ were arbitrary we obtain that $g\in E_{T^n x}(f)$.
So we can conclude that $f(n,x)\cdot E_x(f)\cdot f(n,x)^{-1}\subseteq E_{T^n x}(f)$, and by symmetry follows the assertion.
\end{proof}

\begin{rem}
Obviously a cocycle is recurrent if and only if the identity is in $E(f)$, and the Lemma above shows us that $\mathbf{1}_G\in E_x(f)$ for some $x\in X$ implies that $\mathbf{1}_G\in E_{T^n x}(f)$ for all integers $n$.
As the $T$-orbit of $x$ is dense in $X$, it follows immediately that $\mathbf{1}_G\in E_y(f)$ for all $y\in X$.
Thus $\mathbf{1}_G\in E_x(f)$ for some $x\in X$ implies that $\mathbf{1}_G\in E(f)$ and the cocycle is recurrent.

Another application of the Lemma uses that in an abelian group conjugation does not affect the essential range and thus $E_{T^n x}(f)=E_x(f)$ for all integers $n$.
The fact that every $T$-orbit is dense now implies that $E_x(f)\subseteq E_y(f)$ for every $y\in X$, and by symmetry it follows that $E_x(f)=E_y(f)$ for all $x,y\in X$.
So we obtain that $E_x(f)=E(f)$ for all $x\in X$.
\end{rem}

It is also an easy consequence that $E(f)$ is always a closed \emph{subgroup} of $G$, independently of whether $G$ is abelian or nonabelian, but the situation is much more complicated for $E_x(f)$.
For the study of $E_x(f)$ we need another closely related definition:
\begin{defi}
 For every $x\in X$ we set
\begin{equation*}
P_x(f)=\big\{g:g=\lim_{k\to\infty} f(n_k,x)\enspace\textup{with}\enspace \lim_{k\to\infty}T^{n_k}x=x\enspace\textup{and}\enspace |n_k|\to\infty\big\}.
\end{equation*}
Obviously $P_x(f)$ is contained in $E_x(f)$ and from the continuity of $f(n,\cdot)$ and the cocycle equality (\ref{eq:cocycle2}) it follows that $P_x(f)$ is a closed \emph{sub-semigroup}, supposed that it is nonempty.
It is obvious that $\mathbf{1}_G\in P_x(f)$ if and only if $(x,g)$ is a recurrent point in the skew product for any $g\in G$, i.e. a point which is a limit point of its $\mathbf{T}_f$-orbit, and then it is easily verified that
\begin{equation*}
P_x(f)=\big\{g\in G: (x,g)\in\bar{\textup{O}}_f (x,\mathbf{1}_G)\big\}.
\end{equation*}
\end{defi}

\begin{rems}
(i) The assertion of Lemma \ref{lm:E_orbit} also holds if one replaces $E_x(f)$ by $P_x(f)$ and $E_{T^n x}(f)$ by $P_{T^n x}(f)$ respectively, and the proof is analogous.

(ii) It is an important fact that the assertion (ii) in the Remarks on essential ranges does not hold any more if $E_{x_n}(f)$ is replaced by $P_{x_n}(f)$ and $E_x(f)$ is replaced by $P_x(f)$ respectively.
\end{rems}

\begin{prop}\label{l:E_x P_x}
The set $\mc D(f)= \big\{x\in X: E_x(f)=P_x(f)\big\}$ contains a dense $G_\delta$-set.
If $f$ is recurrent, then for every $x\in\mc D(f)$ the set $P_x(f)=E_x(f)$ is a closed subgroup of $G$.
\end{prop}

\begin{proof}
If $C$ is a compact subset of $G$, then the set $\mc D_C=\big\{x\in X: E_x(f)\cap C=\varnothing\big\}$ is open in $X$.
Indeed, suppose that $E_x(f)\cap C=\varnothing$ while $g_n\in E_{x_n}(f)\cap C$ for a sequence $x_n\to x$.
As $C$ is compact there exists a limit point $g\in C$ of a convergent subsequence $\{g_{n_k}\}_{k\geq 1}$, and a contradiction occurs as this limit point is an element of $E_x(f)\cap C$.

Now let $U$ be a relatively compact open neighbourhood of $\mathbf{1}_G$, and let $g\in G$ and $\varepsilon>0$ be arbitrary.
Let us consider then the following subset of $X$:
\begin{equation*}
\mc G_{(g,U,\varepsilon)}=\mc D_{g\overline{U}}\cup\big\{x: f(n,x)\in gU^2\enspace\textup{for some}\enspace 0\neq n\in\Z\enspace\textup{with}\enspace \delta(x,T^n x)<\varepsilon\big\}
\end{equation*}
This set is open, because $\mc D_{g\overline{U}}$ is open and the second component of the union above is a countable union of open sets.
Furthermore this set is also dense in $X$, because for any $x\in X$ either $x\in\mc D_{g\overline{U}}$ or otherwise, as $g U^2$ is an open neighbourhood of some $h\in E_x(f)$, points out of the second component of the union are arbitrarily close to $x$.
If $\{g_k\}_{k\geq 1}$ is dense sequence in $G$ and $\{U_l\}_{l\geq 1}$ is a neighbourhood base at $\mathbf{1}_G$, then by Baire's theorem the set
\begin{equation*}
\mc G=\bigcap_{k,l,m\geq 1}\mc G_{(g_k,U_l,2^{-m})}
\end{equation*}
is a dense $G_\delta$-subset of $X$.
Let $x\in\mc G$ be arbitrary and fix $g\in E_x(f)$, then there exist increasing sequences of integers $\{k_m\}_{m\geq 1}$ and $\{l_m\}_{m\geq 1}$ so that the sets $g_{k_m} \overline{U}_{l_m}^2$ form a neighbourhood base at $g$.
From $E_x(f)\cap g_{k_m} \overline{U}_{l_m}\neq\varnothing$ we can conclude for any $m\geq 1$ and $\varepsilon>0$ that there is an integer $n\neq 0$ so that $\delta(x,T^n x)<\varepsilon$ and $f(n,x)\in g_{k_m} U_{l_m}^2$.
Therefore we obtain that $g\in P_x(f)$.

For a recurrent cocycle $f$ it follows that $\mathbf{1}_G\in E_y(f)$ for every $y\in X$, and for every $x\in\mc D(f)$ the set $P_x(f)=E_x(f)$ is a nonempty and symmetric sub-semigroup of $G$ and thus it is a subgroup of $G$.
\end{proof}

Now we want to give a definition of regularity for a cocycle, and this is done in analogy to the measure theoretic setting of the problem:
\begin{defi}
A cocycle $f$ is called \emph{regular}, if its skew product transformation $\mathbf{T}_f$ admits an orbit closure of a single point $(x_0,g_0)$ which projects onto all of $X$ under the first projection.
Such an orbit closure will be called a \emph{surjective orbit closure} of $\mathbf{T}_f$.
\footnote{This terminology is influenced by Eli Glasner and Eyal Masad, who use the term \textit{surjective cocycle} for what we call a regular cocycle.}
\end{defi}

We shall later see that for an abelian group $G$ our definition of regularity coincides with the definition given in \cite{LM}, i.e. that the factor cocycle $\tilde f(n,x)= f(n,x)\cdot E(f)$ into $G/E(f)$ does not assume the infinity as an improper essential value.

In the following two sections of this paper we shall investigate the structure of surjective $\mathbf{T}_f$-orbit closures and essential ranges, at first for cocycles with values in arbitrary l.c.s. groups and then with values in nilpotent groups, in which case stronger results can be achieved.
We shall then state a general regularity theorem:
If $T$ is a minimal rotation on a locally connected compact group and $f$ is a recurrent cocycle with values in a nilpotent l.c.s. group, then $f$ is regular.
Together with the results on regular cocycles in the preceding this theorem generalises a result of Atkinson in the paper \cite{A} to nonabelian groups.

According to personal communication Eli Glasner and Eyal Masad developed independently from our work a more abstract approach to the problem of regular cocycles.
In difference to our approach they use a topological analogue of the ergodic decomposition.
In particular, the assertions of Theorem \ref{t:structure}, but except compactness and minimality of $C/H$, as well as results out of Theorem \ref{t:transitive_points} are consequences of their work.

The authors would like to thank Klaus Schmidt for advice and encouragement.

\section{Regular cocycles in locally compact groups}
\label{s:regular}
The study of regular cocycles will be mainly accomplished in the skew product, and an important role plays the continuous action of $G$ on $Y=X\times G$ via the right translations $\{R_h:h\in G\}$ defined by
\begin{equation}
\big(h,(x,g)\big) \mapsto R_h (x,g) = (x,gh^{-1}).
\end{equation}
For any closed subset $C$ of $Y$ we denote the \emph{Stabiliser} with respect to the right translations by
\begin{equation}
\Stab(C)= \big\{g\in G: R_g (C)= C \big\}.
\end{equation}
It is clear from the definition that $\Stab(C)$ is a subgroup of $G$.
Furthermore, if $g_n$ is a sequence in $\Stab(C)$ with $g_n\to g\in G$, then for any $y\in Y$ it follows from the closedness of $C$ that $R_g (y)=\lim_{n\to\infty} R_{g_n} (y)\in C$.
Thus $R_g(C)\subseteq C$, and on the other hand we obtain for any fixed $y\in C$ that $R_{g_n^{-1}}y\to R_{g^{-1}}y$ and hence $R_{g^{-1}}y\in C$.
From $R_g\circ R_{g^{-1}}y=y$ we can conclude that $C\subseteq R_g(C)$, and thus $\Stab(C)$ is a closed subgroup of $G$.

Let $H$ be a closed subgroup of $G$ and denote by $\pi_{G/H}$ the quotient mapping from $G$ onto the homogeneous space $G/H$.
We shall use the notation $\pi_{G/H}$ for the map $(x,g)\mapsto (x,gH)$ from $X\times G$ onto $X\times G/H$.
The skew product transformation $\mathbf{T}_f$ acts on $X\times G/H$ by
\begin{equation}
\mathbf{T}_f (x,gH) = \big( Tx, f(x)\cdot gH \big),
\end{equation}
which shall cause no confusion with $\mathbf{T}_f$ defined on $X\times G$.
In particular we have $\mathbf{T}_f\circ \pi_{G/H} = \pi_{G/H} \circ \mathbf{T}_f$.
Let $C$ be a closed and $\mathbf{T}_f$-invariant subset of $X\times G$ and denote $H=\Stab(C)$, then $C$ contains with any $(x,g)$ the whole left coset $(x,gH)$.
Thus $C$ itself is a subset of $X\times G/H$ and the quotient set $C/H= \pi_{G/H} (C)$ is also closed and $\mathbf{T}_f$-invariant.

\begin{theo}[Structure of surjective orbit closures]\label{t:structure}
Suppose that $f$ is a regular cocycle taking values in a l.c.s. group $G$, and let $C$ be a surjective $\mathbf{T}_f$-orbit closure.
Then for every point $y=(x,g)\in C$ with a dense orbit in $C$ it follows that
\begin{equation}\label{e:single_coset}
C_x=\big\{g'\in G: (x,g')\in C\big\} = g H,
\end{equation}
in which $H=\Stab(C)$.
Furthermore, the quotient set $C/H = \pi_{G/H} (C)$ is compact and $\mathbf{T}_f:C/H\longrightarrow C/H$ is a minimal homeomorphism.
In particular, for every $x\in X$ there exists a compact set $K_x$ in $G$ so that $C_x= K_x H$.
\end{theo}

\begin{rem}
As a consequence of the following Theorem \ref{t:transitive_points} the equation (\ref{e:single_coset}) holds on the set $\mc D(f)=\big\{x\in X: E_x(f)=P_x(f)\big\}$, which contains a dense $G_\delta$ set.
Later we shall prove for nilpotent groups $G$ that $C_x=g_x H$ for every $x\in X$.
\end{rem}

For the proof of the statement above the following Lemmas are necessary.

\begin{lem}\label{l:pseudo_open}
Suppose that $C$ is surjective $\mathbf{T}_f$-orbit closure of a regular cocycle $f$.
Then every nonempty, relatively open subset $\mc O_C$ of $C$ projects onto a set with nonempty interior under the first projection $\pi_X$.
\end{lem}

\begin{proof}
We show first that for any open neighbourhood $U$ of $\mathbf{1}_G$ there is an open $\mathbf{T}_f$-invariant subset $\mc O\subseteq R_U(C)=\cup_{g\in U}R_g(C)$ so that $\mc O\cap C$ is nonempty and dense in $C$.
Let $V$ be any relatively compact open neighbourhood of $\mathbf{1}_G$ so that $\overline V^{-1} \overline V\subseteq U$ and choose a sequence $\{h_n\}_{n\geq 1}$ which is dense in $G$.
Obviously it follows that $\bigcup_{n\geq 1} R_{h_n \overline V} C = X\times G$, and as $\overline V$ is compact and the right translations are isometries each of the sets $R_{h_n \overline V} C$ is closed.
Hence by Baire's theorem there is an integer $n$ so that $R_{h_n \overline V} C$ as well as $R_{\overline V} C = R_{h_n^{-1}} \big(R_{h_n \overline V} C\big)$ has a nonempty interior.
Thus $R_{\overline V} C$ contains a $\mathbf{T}_f$-invariant open set $\mc O'$.
Note that $\mc O'\cap C$ might be empty, however $R_{\overline V^{-1} \,\overline V} C$ contains the open set $\mc O=R_{\overline V^{-1}} \mc O'$ which has nonempty intersection with $C$.

Now suppose that $\mc O$ is open in $X\times G$ and that $\mc O\cap C\neq\varnothing$.
Choose a smaller open set $\mc O'$ and an open neighbourhood $U$ of $\mathbf{1}_G$ so that $\mc O'\cap C\neq \varnothing$ and $R_{U^{-1}} \mc O'\subseteq \mc O$.
By the preceding argument $R_U C$ contains an open set which is relatively dense in $C$, and thus $\big(R_U C \big)^\circ \cap \mc O' \neq \varnothing$.
But this implies that the set $R_{U^{-1}}\mc O' \cap C$ is open in $C$ and its projection has a nonempty interior.
\end{proof}

\begin{lem}\label{l:nice_point}
Let $C$ be a surjective $\mathbf{T}_f$-orbit closure of a regular cocycle $f$.
Then there exists a point $(x_0,g_0)\in C$ with a dense orbit in $C$ so that $E_{x_0}(f)=P_{x_0}(f)$.
\end{lem}

\begin{proof}
By Lemma \ref{l:E_x P_x} the set $\mc D(f)=\big\{x\in X: E_x(f)=P_x(f)\big\}$ contains a dense $G_\delta$ set $D$.
It is clear that $D_1=\pi_X^{-1} (D)\cap C$ is also a $G_\delta$ set and Lemma \ref{l:pseudo_open} implies that $D_1$ is also dense in $C$.
As $C$ is topologically transitive, the set $D_2$ of all topologically transitive points forms a dense $G_\delta$-set in $C$ (cf. \cite{G-H}, Theorem 9.20).
Then by Baire's theorem $D_1\cap D_2$ is nonempty and hence we can find a point $y_0=(x_0,g_0)\in C$ for which both $\bar{\textup{O}}_f(y_0) = C$ and $x_0\in\mc D(f)$ hold.
\end{proof}

\begin{lem}\label{l:compact}
Suppose that $C_H\subseteq X\times C/H$ is a closed $\mathbf{T}_f$-invariant set which projects onto all of $X$.
Then there exists a compact set $K\subseteq G$ with the property that
\begin{equation*}
X= \pi_X \big(C_H \cap (X\times KH) \big).
\end{equation*}
\end{lem}

\begin{proof}
Let $\{K_n\}_{n\geq 1}$ be a sequence of compact sets for which $G/H=\bigcup_n K_n H$.
Then all the sets $A_n = \pi_X \big( C_H \cap (X\times K_n H) \big)$ are compact subsets of $X$ and $X=\bigcup_n A_n$.
By Baire's theorem there exists an open set $\mc O\subseteq X$ which is contained in some $A_{n_0}$.
As $C_H$ is $\mathbf{T}_f$-invariant, it follows for all integers $k$ that
\begin{equation*}
T^k (\mc O) \subseteq \pi_X \big ( C_H \cap \mathbf{T}_f^k (X \times K_{n_0} H)\big),
\end{equation*}
but as $X$ is compact and minimal we can find a positive integer $N$ so that $X=\bigcup_{k=1}^N T^k (\mc O)$.
Thus the compact set $K= \bigcup_{k=1}^N f(k,X)\cdot K_n$ is sufficient for the equality that $\pi_X\big(C_H \cap (X\times KH)\big)=X$.
\end{proof}

\begin{proof}[Proof of Theorem \ref{t:structure}]
Lemma \ref{l:nice_point} shows that there is a point $(x_0,g_0)\in C$ with a dense $\mathbf{T}_f$-orbit in $C$ and $E_{x_0}(f)=P_{x_0}(f)$.
In particular, $P_{x_0}(f)$ is a subgroup of $G$ and the vertical section
\begin{equation*}
C_{x_0}= P_{x_0}(f)\cdot g_0 = g_0 g_0^{-1}\cdot P_{x_0}(f) \cdot g_0
\end{equation*}
is a left coset of the subgroup $H'=g_0^{-1}\cdot P_{x_0}(f)\cdot g_0$.
It follows from the definition of the stabiliser that $\Stab(C)\subseteq H'$, and in order to show that $H'=\Stab(C)$ we choose an arbitrary $h\in H'$ and observe that
\begin{equation*}
R_h (C) = R_h \bar{\textup{O}}_f(x_0,g_0) = \bar{\textup{O}}_f \big(R_h (x_0,g_0)\big) \subseteq C,
\end{equation*}
since $R_h (x_0,g_0)\in C$.
As $H'$ is a subgroup of $G$, we have $R_h(C)=C$ for every $h\in H'$ and thus $H'\subseteq \Stab(C)$.

From Lemma \ref{l:compact} it follows that there exists a compact set $K\subset G$ so that $C/H \cap (X\times KH)$ projects onto all of $X$.
But for any integer $n$ the vertical section of $C/H$ at $T^n x_0$ consists of one left coset of $H$, which therefore must be contained in $K H$.
Thus the orbit of $(x_0,g_0 H)$ under $\mathbf{T}_f$ is a subset of $X\times K H$, and as this orbit is also dense in $C/H$ it follows that $C/H\subseteq X\times K H$.

For the proof of the minimality of $\mathbf{T}_f$ on $C/H$ we fix a point $(x,gH)\in C/H$ and choose a sequence of integers $\{n_k\}_{k\geq 1}$ for which $T^{n_k}x \rightarrow x_0$.
As $C/H$ is compact, a subsequence of $\mathbf{T}_f^{n_k} (x,gH)$ converges to a point $y\in C/H$ with $\pi_X (y) = x_0$.
But the only point in $C/H$ which projects onto $x_0$ is $(x_0,g_0H)$, and the fact that the orbit of $(x_0,g_0H)$ under $\mathbf{T}_f$ is dense in $C/H$ implies that $\mathbf{T}_f$ is minimal.

It remains to prove that for any $(x,g)\in C$ with a dense orbit in $C$ the vertical section $C_x$ consists only of a single left coset.
As the orbit of $(x,g)$ is dense in $C$, we can find an integer sequence $\{n_k\}_{k\geq 1}$ with $T^{n_k}x\rightarrow x_0$ and $f(n_k,x)\cdot g\rightarrow g_0$.
Thus it follows for any $h\in C_x$ that $f(n_k,x)\cdot h \rightarrow g_0 g^{-1} h$ and we can conclude that $g_0 g^{-1} C_x \subseteq C_{x_0}$.
By symmetry the same inclusion holds if we interchange $(x,g)$ and $(x_0,g_0)$, and ultimately we obtain that $C_x= g  g_0^{-1} C_{x_0}= g H$.
\end{proof}

The following theorem shows that the topologically transitive points of $C$ can be characterised by vertical sections and essential ranges.

\begin{theo}[Essential ranges and transitive points]\label{t:transitive_points}
Suppose that $f$ is a regular cocycle taking values in a l.c.s. group $G$ and let $C$ be a surjective $\mathbf{T}_f$-orbit closure.
For any $x\in X$ we have the inclusions
\begin{equation}\label{e:C_x_er}
P_x(f)\subseteq C_x C_x^{-1}=\{h k^{-1}:h,k\in C_x\}\subseteq E_x(f),
\end{equation}
and equality between these three sets holds if and only if the orbit of some $y\in C$ with $\pi_X(y)=x$ is dense in $C$.
Furthermore, for any $(x,g)\in C$ with a dense orbit in $C$ we have the equality that
\begin{equation}
E_x(f)=g H g^{-1},
\end{equation}
in which $H=\Stab(C)$.
Thus for all $x$ in the set $\mc D(f)=\big\{x\in X: E_x(f)=P_x(f)\big\}$, which contains a dense $G_\delta$ set, the essential range $E_x(f)$ is conjugate to the closed subgroup $H$ and $\bar{\textup{O}}_f(x,\mathbf{1}_G)$ is a surjective  $\mathbf{T}_f$-orbit closure.
\end{theo}

\begin{proof}
The inclusion $P_x(f)\subseteq C_x C_x^{-1}$ follows immediately from $P_x(f)\cdot g\subseteq C_x$ for any $g\in C_x$.
Assume that $h,k\in C_x$ and let $(x_0,g_0)\in C$ be a point with a dense orbit in $C$.
Then we can find two sequences of integers $\{m_k\}_{k\geq 1}$ and $\{n_k\}_{k\geq 1}$ so that $\mathbf{T}_f^{m_k} (x_0,g_0)\to (x,h)$ and $\mathbf{T}_f^{n_k} (x_0,g_0)\to (x,k)$.
But this means that both $T^{m_k} x_0\to x$ and $T^{n_k} x_0\to x$ while $f(n_k-m_k,T^{m_k} x_0)\to k h^{-1}$, and hence we can conclude that $\{k h^{-1}: h,k \in C_x\}\subseteq E_x(f)$.

If $(x,g)\in C$ is point with $E_x(f)=P_x(f)$, then it follows from $\bar{\textup{O}}_f(x,g)\subseteq C$ that $P_x(f)\subseteq C_x g^{-1}\subseteq E_x(f)$, and hence $P_x(f)=C_x g^{-1}$.
From $g\in C_x$ we obtain the inclusion $gH\subseteq P_x(f)\cdot g=C_x$, and from the fact that $\mathbf{T}_f$ is minimal on $C/H$ we can conclude that $\bar{\textup{O}}_f(x,g)=C$.

We assume now that $\bar{\textup{O}}_f(x,g)= C$ and choose according to Lemma \ref{l:nice_point} a point $(x_0,g_0)\in C$ with a dense orbit in $C$ and $E_{x_0}(f)=P_{x_0}(f)$.
We let $\{n_k\}_{k\geq 1}$ be a sequence with ${\mathbf{T}_f^{n_k}(x,g)\to (x_0,g_0)}$ and apply Lemma \ref{lm:E_orbit} to conclude that $E_{T^{n_k}x}(f) = f(n_k,x)\cdot E_{x}(f)\cdot f(n_k,x)^{-1}$.
But as $f(n_k,x)\to g_0 g^{-1}$ while $T^{n_k}x\to x_0$ it follows from the definition of the essential range that $g_0 g^{-1}\cdot E_{x}(f)\cdot g g_0^{-1}\subseteq E_{x_0}(f)$.
At this point we only used that the orbit of $(x,g)$ is dense in $C$ and as the orbit of $(x_0,g_0)$ is also dense in $C$, we obtain the converse inclusion by symmetry and thus
\begin{equation*}
E_x(f) =  g g_0^{-1}\cdot E_{x_0}(f)\cdot g_0 g^{-1}.
\end{equation*}
From $f(n_k,x)\to g_0 g^{-1}$ and $T^{n_k}x\to x_0$ it follows that $g_0 g^{-1}C_x\subseteq C_{x_0}$, and again by symmetry we can conclude that
\begin{equation*}
C_x = g g_0^{-1} C_{x_0} .
\end{equation*}
As both $(x,g)$ and $(x_0,g_0)$ have a dense orbit in $C$, we obtain that $C_x=P_x(f)\cdot g$ and $C_{x_0}=P_{x_0}(f)\cdot g_0$, and together with the equality $E_{x_0}(f)=P_{x_0}(f)$ we can conclude that $C_{x_0}=P_{x_0}(f)\cdot g_0=E_{x_0}(f)\cdot g_0$.
Thus we have that
\begin{equation*}
P_x(f)=C_x g^{-1} = g g_0^{-1}\cdot E_{x_0}(f)\cdot g_0 g^{-1}=E_x(f).
\end{equation*}

Furthermore, if the orbit of $(x,g)$ is dense in $C$ then Theorem \ref{t:structure} states that $C_x=P_x(f)\cdot g=gH$, and hence $E_x(f)=P_x(f) = g H g^{-1}$.
\end{proof}

Let us now discuss some consequences of the Theorems \ref{t:structure} and \ref{t:transitive_points}.
If $f$ is a regular cocycle and $x_0$ is any point out of $\mc D(f)$, then the $\mathbf{T}_f$-orbit closure of $(x_0,\textbf{1}_G)$ is compact modulo the subgroup $E_{x_0}(f)$.
Conversely, if $x_0\in \mc D(f)$ is a point so that $C=\bar{\textup{O}}_f(x_0,\mathbf{1}_G)$ is compact modulo $E_{x_0}(f)$, then its projection $\pi_X(C)$ is a nonempty compact subset of $X$ and from the minimality it follows that $\pi_X(C)=X$.
Thus we can state the following characterisation of regularity:
\begin{cor}\label{cor:comp}
A cocycle $f$ is regular, if and only if for one (and therefore every) point $x_0$ belonging to $\mc D(f)=\big\{x\in X: E_x(f)=P_x(f)\big\}$ the $\mathbf{T}_f$-orbit closure of $\big(x_0,E_{x_0}(f)\big)$ in $X\times G/E_{x_0}(f)$ is compact.
\end{cor}
For an abelian group $G$ we already know that $E_x(f)=E(f)$ for every $x\in X$.
Thus the cocycle $f$ is regular if and only its factor cocycle $\tilde f(n,\cdot)= f(n,\cdot)\cdot E(f)$ into $G/E(f)$ does not assume the infinity as an extended essential value.
This property is obviously equivalent to the compactness of the $\mathbf{T}_{\tilde f}$-orbit closure of any point $(x_0,E(f))$, and thus the notion of regularity given in this paper generalises the notion used in \cite{LM}.

Every surjective $\mathbf{T}_f$-orbit closure is a right translate of any other.
In fact, if $C_1$ and $C_2$ are surjective $\mathbf{T}_f$-orbit closures then for any point $x\in\mc D(f)$ we have $C_i=\bar{\textup{O}}_f(x,g_i)$ for appropriate elements $g_i\in G$.
Therefore $C_2= R_{g_2^{-1}g_1} C_1$ is a right translate of $C_1$ and their associated stabiliser subgroups $H_i=\Stab(C_i)$ are always conjugate.

It is an important question whether $X\times G$ is a disjoint union of surjective $\mathbf{T}_f$-orbit closures.
This is the case if and only if \emph{every} vertical section $C_x$ of an orbit closure $C$ consists only of one left coset of its stability group $H$.
In the next section it will be shown that this condition is always fulfilled if $G$ is nilpotent, but it is unclear to the authors whether this is also true for an arbitrary l.c.s. group.

\section{Regular cocycles in nilpotent locally compact groups}
\label{s:regular_nil}

Suppose that $G$ is a nilpotent l.c.s. group and denote by $Z(G)$ the centre of $G$, which is a closed normal abelian subgroup of $G$.
If we define a sequence of l.c.s. groups inductively by $G_0=G$ and $G_{n+1}=G_n/Z(G_n)$, then $G_n$ is trivial for some positive integer $n$.
We say that $G$ is $n$-step nilpotent if $n$ is the smallest integer for which $G_n$ is trivial.
It is easy to see that closed subgroups and homomorphic images of nilpotent groups are also nilpotent.
In the main results of the paper we shall use the following well known properties of nilpotent groups:

\begin{lem}\label{l:double_cosets}
Let $G$ be a nilpotent l.c.s. group and let $H$ be a closed subgroup of $G$.
Then for every $g\notin H$ the closure of the double coset $HgH$ does not contain the identity $\mathbf{1}_G$.
Furthermore, for any $g\neq\mathbf{1}_G$ the closure of the conjugation class $C(g)=\{hgh^{-1}:h\in G\}$ of $g$ does not contain $\mathbf{1}_G$.
\end{lem}

\begin{proof}
Our statement is obviously true for any $1$-step nilpotent (abelian) group $G$.
Suppose that the statement of our lemma is true for all $n$-step nilpotent groups and let $G$ be $(n+1)$-step nilpotent with a closed subgroup $H$.
We denote by the projection of $G$ onto $\tilde{G}=G/Z(G)$ by $\pi$.

If we have $\pi(g) \in\tilde{H}=\overline{\pi(H)}$, then for every open neighbourhood $U$ of $\mathbf{1}_G$ we can find $h\in H$ so that $h\in U g\cdot Z(G)$, because $\pi$ is an open mapping.
It follows that there exists an element $k\in Z(G)$ with $k H\subseteq U g H$, and as $k$ is in the centre of $G$ any left translation by an element of $H$ leaves the set $kH$ invariant.
Hence $gH$ is the limit point of a sequence of fixed points under the left translations by elements of $H$, and thus $gH$ is also a fixed point.
We can conclude that $\overline{HgH}=HgH= gH$, and this set does not contain $\mathbf{1}_G$ if $g\notin H$.

Otherwise $\pi(g)\notin\tilde{H}$ and as $\tilde{G}$ is $n$-step nilpotent the identity in $\tilde{G}$ is not contained in the closure of $\tilde{H}\pi(g)\tilde{H}$, and from the continuity of $\pi$ it follows that also $\mathbf{1}_G \notin \overline{H g H}$.

The proof of the second statement is similar: If $g\in Z(G)$ then $C(g)=\{g\}$, otherwise $\pi(g)\neq\mathbf{1}_{\tilde{G}}$ and the statement follows from the induction hypothesis on $\tilde{G}$ and the continuity of $\pi$.
\end{proof}

\begin{theo}[Structure of surjective orbit closures]\label{t:structure_nil}
Suppose that $f$ is a regular cocycle taking its values in a nilpotent l.c.s. group $G$, and let $C$ be a surjective $\mathbf{T}_f$-orbit closure.
We have then the following statements, in which $H$ denotes the closed subgroup $\Stab(C)$:
\begin{enumerate}
\item For every $x\in X$ there is some $g_x\in G$ so that $C_x=g_x H$.
\item The mapping $\gamma:X\longrightarrow G/H$ with $\gamma(x)= g_x H$ is continuous.
\end{enumerate}
\end{theo}

\begin{proof}
Suppose that $C=\bar{\textup{O}}_f(x_0,g_0)$ is a surjective orbit closure, let $x\in X$ be an arbitrary point and let $\{n_k\}_{k\geq 1}$ be a sequence with $T^{n_k} x\rightarrow x_0$.
From Theorem \ref{t:structure} we know that $C/H$ is compact and that $C_{x_0}=g_0 H$.
If we choose a neighbourhood base $\{U_k\}_{k\geq 1}$ at $\mathbf{1}_G$, then there is a subsequence $\{m_k\}_{k\geq 1}$ of $\{n_k\}_{k\geq 1}$ so that
\begin{equation*}
f(m_k,x)\cdot C_x =C_{T^{m_k}x} \subseteq U_k g_0 H \enspace\textup{for all}\enspace k\geq 1.
\end{equation*}
Indeed, assume that there exists a neighbourhood $U_l\in\{U_k\}_{k\geq 1}$ and a sequence $\{g_k H\}_{k\geq 1}\subseteq C_x$ so that $f(n_k,x)\cdot g_k H \notin U_l g_0 H$ for infinitely many $k\geq 1$.
As $\mathbf{T}_f$ is a minimal homeomorphism on the compact space $C/H$ and $T^{n_k} x\to x_0$, it follows that there exits a least one limit point of the sequence $f(n_k,x)\cdot g_k H$ outside of $U_l g_0 H$.
But such a limit point is an element of $C_{x_0}/H$ apart from $g_0 H$, and this leads to a contradiction.

If we assume that $g,g' \in C_x$ then
\begin{equation*}
g^{-1} g'=g^{-1}\cdot f(m_k,x)^{-1}\cdot f(m_k,x) \cdot g'\in H g_0^{-1} U_k^{-1} U_k g_0 H ,
\end{equation*}
and hence for suitably chosen sequences $\{h_k\}_{k\geq 1}, \{h_k'\}_{k\geq 1}\subseteq H$ it follows that
\begin{equation*}
h_k^{-1} g^{-1} g' h_k'\in g_0^{-1} U_k^{-1} U_k g_0=V_k.
\end{equation*}
The open sets $\{V_k\}_{k\geq 1}$ also define a neighbourhood base at $\mathbf{1}_G$ and thus $\mathbf{1}_G\in \overline{H g^{-1} g' H}$, and then Lemma \ref{l:double_cosets} implies that $g^{-1} g'\in H$.
But as $g, g'\in C_x$ were arbitrary the set $C_x$ consists of just one left $H$-coset.

The projection $\pi_X$ from $C/H$ to $X$ is a continuous, open, one to one, and onto mapping, and hence it is a homeomorphism of the compact spaces $C/H$ and $X$.
As $\pi^{-1}_X$ is continuous the proof is complete.
\end{proof}

\begin{rem}
Theorem \ref{t:structure_nil} as well as all other results of this section actually hold for any group $G$ satisfying the double coset property in Lemma \ref{l:double_cosets}.
For instance, every discrete group or every group which admits a biinvariant metric obviously satisfies this property.
Note that the double coset property is weaker than nilpotency.
For example, if the sequence of groups $G_n$ defined by $G_0=G$ and $G_{n+1}=G_n/Z(G_n)$ terminates with a group which admits a biinvariant metric (e.g. any compact group), then $G$ also satisfies the double coset property.
The proof is then analogous to the proof of Lemma \ref{l:double_cosets}.
\end{rem}

Let $H$ be any closed subgroup of $G$, and denote the set of all subgroups which are conjugate to $H$ over an element of $G$ by $H^G=\{gHg^{-1}:g\in G\}$.
If we let $\Stab(H)=\{g\in G: gHg^{-1}= H\}$ be the stabiliser of $H$ in $G$, then the mapping
\begin{equation}
\varphi:G/\Stab(H)\longrightarrow H^G\enspace{with}\enspace \varphi(g \cdot\Stab(H))= gHg^{-1}
\end{equation}
is a bijection between $G/\Stab(H)$ and $H^G$.
We shall identify $H^G$ with $G/\Stab(H)$ and turn it into a locally compact topological space by means of this identification.

\begin{theo}[Essential ranges]\label{t:essential_nil}
Suppose that $f$ is a regular cocycle taking its values in a nilpotent l.c.s. group $G$, let $C$ be a surjective $\mathbf{T}_f$-orbit closure and let $H=\Stab(C)$.
For every $x\in X$ we have the equality that
\begin{equation}
E_x(f)=C_x C_x^{-1}=g_x H g_x^{-1},
\end{equation}
in which $g_x$ is determined as in Theorem \ref{t:structure_nil}.
It follows that $E_x(f)\in H^G$ for all $x\in X$ and that $x\mapsto E_x(f)$ is a continuous map from $X$ into $H^G$.
\end{theo}

\begin{proof}
We want to show first that $\pi_X$ is an open mapping from $C$ onto $X$, and from Theorem \ref{t:structure_nil} we already know that the first projection is an open mapping from $C/H$ onto $X$.
As every vertical section of $C$ consists of only one left coset of $H$ we can conclude for every open subset $\mc O\subseteq X\times G$ with $\mc O\cap C\neq\varnothing$ that
\begin{equation*}
\pi_{G/H} (\mc O\cap C) = \pi_{G/H}(\mc O)\cap C/H ,
\end{equation*}
where the latter set is open in $C/H$.
Therefore $\pi_X(\mc O\cap C) = \pi_X \big(\pi_{G/H} (\mc O\cap C)\big)$ is open in $X$.

The inclusion $C_x C_x^{-1}\subseteq E_x(f)$ is part of Theorem \ref{t:transitive_points}, and to prove the converse inclusion let $(x,g)\in C$ and $h\in E_x(f)$ be arbitrary.
Let $U$ be an arbitrary open neighbourhood of $\mathbf{1}_G$ and choose an open neighbourhood $V$ of $\mathbf{1}_G$ with $VhVh^{-1}\subseteq U$.
For any $\delta>0$ the set $\pi_X\big(C\cap (B(x,\delta)\times V g)\big)$ is open and contains an open ball $B(x,\delta')$ for some $\delta'>0$.
As $h\in E_x(f)$, we can find a point $x'$ contained in $B(x,\delta')$ so that $T^n x'\in B(x,\delta')$ and $f(n,x')\in Vh$.
If $(x',g')$ is a point in $C\cap (B(x,\delta)\times V g)$ which projects on $x'$, then
\begin{equation*}
\mathbf{T}_f^n (x',g') = \big(T^n x', f(n,x')g'\big ) \in B(x,\delta')\times Vhg' \subseteq B(x,\delta)\times U hg.
\end{equation*}
As $C$ is closed while $U$ and $\delta$ were arbitrary we can conclude that $(x,hg)\in C$ and thus $h\in C_x C_x^{-1}$.

Theorem \ref{t:structure_nil} says that $C_x= g_x H$ and that the map $\gamma:X\longrightarrow G/H$ with $\gamma(x)= g_x H$ is continuous, and therefore $E_x(f)=C_x C_x^{-1}$ is equal to the conjugate group $g_x H g_x^{-1}\in H^G$, which corresponds to the element $g_x\cdot\Stab(H)$ in $G/\Stab(H)$.
The projection $p:G/H\longrightarrow G/\Stab(H)$ with $p(gH)= g\cdot\Stab(H)$ is continuous and thus the map $x\mapsto E_x(f) = \varphi\circ p \circ\gamma (x)$ is also continuous, in which $\varphi$ denotes the homeomorphism between $G/\Stab(H)$ and $H^G$.
\end{proof}

\section{Regularity of cocycles for minimal rotations}
We suppose from now on that $X$ is a locally connected compact group with a minimal rotation $T$, and we let $\delta(\cdot,\cdot)$ denote an invariant metric on $X$.
It follows that $T$ is a minimal isometry with respect to $\delta$, and $\delta(T^k x,x)<\varepsilon$ for some $x\in X$ implies that  $\delta(T^k y,y)<\varepsilon$ for all $y\in X$.
Thus integer sequences $\{k_n\}_{n\geq 1}$ with $T^{k_n}x\to x$ for some $x\in X$ fulfil that $T^{k_n}\to \textup{Id}_X$ uniformly, and such sequences are often called rigidity times.

\begin{defi}
Let $H$ be a closed subgroup of $G$ and let $H^G$ denote the set of all subgroups conjugate to $H$.
We say that a continuous map $y\mapsto H_y$ from $X$ into $H^G$ is a \emph{consistent selection} of subgroups in the essential ranges of the cocycle $f$, if it fulfils that
\begin{equation}\label{eq:H_E}
H_x\subseteq E_x(f)
\end{equation}
for all $x\in X$ and in analogy to the essential ranges that
\begin{equation}\label{eq:H_ORBIT}
H_{T^n x}=f(n,x)\cdot H_x \cdot f(n,x)^{-1}
\end{equation}
for all $x\in X$ and all integers $n$.
\end{defi}

In the following two lemmas we prove some elementary properties of such consistent selections of subgroups.

\begin{lem}\label{lem:eo}
Let $U$ be a relatively compact open neighbourhood of $\mathbf{1}_G$ so that $E_z(f)\cap (\bar{U} H_z\smallsetminus U H_z)=\varnothing$ for some $z\in X$.
Then there exists an $\varepsilon >0$ so that for all $y\in X$ and $n\in\mathbb{Z}$ with $\delta (y,z)<\varepsilon$ and  $\delta (T^n y,z)<\varepsilon$ it holds that
\begin{equation}\label{eq:eo}
\textup{either}\enspace f(n,y) \cdot H_y\cap U H_y\neq\varnothing\enspace\textup{or}\enspace f(n,y) \cdot H_y\cap
\bar{U}H_y=\varnothing.
\end{equation}
Furthermore, for every $y\in X$ the set $E_y(f)$ is a subset of the quotient space $G/H_y$, i.e for any $g\in E_y(f)$ we have $g H_y\subseteq E_y(f)$.
\end{lem}

\begin{proof}
Assume that there exists a sequence $\{(n_k,y_k)\}_{k\geq 1}\subseteq\mathbb{Z}\times X$ with $y_k\to z$ and $T^{n_k}y_k\to z$ as well as $f(n_k,y_k) \cdot H_{y_k}\cap U=\varnothing$ and $f(n_k,y_k) \cdot H_{y_k}\cap\bar{U} \neq\varnothing$ for all integers $k\geq 1$.
We select then a sequence $g_k\in f(n_k,y_k) \cdot H_{y_k}\cap(\bar{U} \smallsetminus U)$ and assume by the compactness of $\bar{U} \smallsetminus U$ that $g_k$ is convergent to $g'\in \bar{U} \smallsetminus U$.
If this limit point $g'$ were within $U H_y$, then $g' h\in U$ for some $h\in H_y$ would imply that $V W \subseteq U$ for an open neighbourhood $V$ of $g'$ and an open neighbourhood $W$ of $h$.
As $g_k\in V$ and $W\cap H_{y_k}\neq\varnothing$ for all $k$ large enough a contradiction to $f(n_k,y_k) \cdot H_{y_k}\cap U=\varnothing$ would occur, and thus $g'\notin U H_y$.
The inclusion $H_{y_k}\subseteq E_{y_k}(f)$ holds for every integer $k\geq 1$ and hence for every fixed $k$ we can find a sequence $\{(m^k_l,y^k_l)\}_{l\geq 1}\subseteq X\times \mathbb{Z}$ so that $y^k_l\to y_k$, $T^{m^k_l}y^k_l\to y_k$ and $f(n_k,y_k) \cdot f(m^k_l,y^k_l)\to g_k$ as $l\to\infty$.
By the continuity of $f(n_k,\cdot)$ we can state that $f(n_k,T^{m^k_l}y^k_l) \cdot f(m^k_l,y^k_l)=f(n_k+m^k_l,y^k_l)\to g_k$ as $l\to\infty$, and in the limit $k\to\infty$ it follows that $g'\in E_z(f)\cap (\bar{U} \smallsetminus U H_y)$, which contradicts the assumptions of the Lemma.

The second assertion of the lemma also follows from the preceding argument, if we let $f(n_k,y_k)\to g\in E_y(f)$ and, as the map $y\mapsto H_y$ is continuous, choose a sequence $\{h_k\}_{k\geq 1}$, $h_k\in H_{y_k}$ with $h_k\to h\in H_y$ and let $f(m^k_l,y^k_l)\to h_k$ as $l\to\infty$.
\end{proof}

\begin{lem}\label{lm:cont}
Let $U\subseteq G$ be an open subset and $C\subseteq G$ a compact subset.
Then for any fixed integer $n$ the sets $\{y\in X: f(n,y) \cdot H_y\cap U H_y\neq\varnothing\}$ and $\{y\in X: f(n,y) \cdot H_y\cap C H_y=\varnothing\}$ are both open.
\end{lem}

\begin{proof}
If $f(n,y)\cdot h\in U$ with some $h\in H_y$, then we can choose two suitable open neighbourhoods $V$ and $W$ of $f(n,y)$ and $h$ respectively so that $VW\subseteq U$.
From the continuity of the maps $y\mapsto f(n,y)$ and $y\mapsto H_y$ it follows that there exists an open neighbourhood $\mathcal{W}$ of $y$ with $f(n,y')\in V$ and $H_{y'}\cap W\neq\varnothing$ for all $y'\in\mathcal{W}$, and thus $f(n,y') \cdot H_{y'}\cap U\neq\varnothing$ for all $y'\in\mathcal{W}$.
For the proof of the second assertion we assume that $f(n,y) \cdot H_y\cap C=\varnothing$, while there is a sequence $\{y_k\}_{k\geq 1}$ convergent to $y$ so that $f(n,y_k) \cdot H_{y_k}\cap C\neq\varnothing$.
For every positive integer $k$ we select then $h_k\in H_{y_k}$ and $g_k=f(n,y_k) \cdot h_k\in C$, and by the compactness of $C$ we can assume that $g_k\to g'\in C$.
But as $f(n,\cdot)$ is continuous the sequence $\{h_k\}_{k\geq 1}$ defined by $h_k=f(n,y_k)^{-1}\cdot g_k$ converges to $f(n,y)^{-1}\cdot g'$.
As the map $y\mapsto H_y$ is continuous and $y_k\to y$ there is a sequence $\{(h'_k,b_k)\}_{k\geq 1}\subseteq H_z\times G$ with $b_k\to\mathbf{1}_G$ so that $h_k=b_k h'_k b_k^{-1}$, and hence $h'_k$ is also convergent to $f(n,y)^{-1}\cdot g$.
We obtain that $f(n,y)^{-1}\cdot g\in H_y$ and a contradiction occurs.
\end{proof}

The following proposition is a generalisation of a result in the paper \cite{A}, and the fact that $X$ is locally connected will be essential.

\begin{prop}\label{pr:REG}
Suppose that $T$ is a minimal rotation on a locally connected compact group $X$ and that $f:X\longrightarrow G$ is a continuous recurrent cocycle for $T$ with values in a l.c.s. group $G$.
Let $y\mapsto H_y$ be a consistent selection of subgroups from $X$ into $H^G$ and assume that for every $y\in X$ there is a neighbourhood base $\{U_{n,y}\}_{n\geq 1}$ of $\mathbf{1}_G$, which consists of relatively compact open sets, so that
\begin{equation*}
E_y(f)\cap (\bar{U}_{n,y} H_y\smallsetminus U_{n,y} H_y)=\varnothing
\end{equation*}
for all $n\geq 1$.
Then for every $x\in\mathcal{D}(f)$ the $\mathbf{T}_f$-orbit closure of $(x,H_x)$ in $G/H_x$ is compact.
\end{prop}

\begin{proof}
We suppose that $(z, g H_x)$ is a point in the $\mathbf{T}_f$-orbit closure of $(x,H_x)$, and as $x\in\mathcal{D}(f)$ there exists an integer sequence $\{n_k\}_{k\geq 1}$ with $|n_k|\rightarrow \infty$ so that $T^{n_k}x\to z$ and $f(n_k,x)\to g$ as $k\to\infty$.
Furthermore it follows from $f(n_k,x)\to g$ and equation (\ref{eq:H_ORBIT}) in the limit that  $H_z=g H_x g^{-1}$.
We fix $\rho>0$ and a positive integer $l$ and apply Lemma \ref{lem:eo} to find
a positive $\varepsilon\leq\rho$ so that for all $y\in X$ and $m\in\mathbb{Z}$ with $\delta (y,z)<\varepsilon$ and $\delta (T^m y,z)<\varepsilon$ we have that
\begin{equation*}
\textup{either}\enspace f(m,y) \cdot H_y\cap U_{l,z}\neq\varnothing\enspace\textup{or}\enspace f(m,y) \cdot H_y\cap \bar{U}_{l,z}=\varnothing.
\end{equation*}
There exists an open connected neighbourhood $\mathcal{O}$ of $z$ with $\mathcal{O}\subseteq\{y:\delta (y,z)<\varepsilon/2\}$, because $X$ is locally connected, and we choose an integer $k'\geq 1$ so that for all $k\geq k'$ we have that $T^{n_{k}}x\in\mathcal{O}$ and
\begin{align*}
f(n_k,x)\cdot f(n_{k'},x)^{-1} =f(n_k-n_{k'},T^{n_{k'}}x)\in U_{l,z},\\
f(n_{k'},x)\cdot f(n_k,x)^{-1} = f(n_{k'}-n_k,T^{n_k}x)\in U_{l,z}.
\end{align*}
The rotation $T$ is an isometry and thus $\delta(T^{n_{k'}-n_k}y,y)<\varepsilon$ and $\delta(T^{n_k-n_{k'}}y,y)<\varepsilon$ uniformly for all $y\in X$, and from the connectedness of $\mathcal{O}$ and Lemma \ref{lm:cont} we can conclude for all $k\geq k'$ that
\begin{equation*}
f(n_{k}-n_{k'},z) \cdot H_z \cap U_{l,z} \neq\varnothing\enspace\textup{and}\enspace f(n_{k'}-n_{k},z) \cdot H_z \cap U_{l,z} \neq\varnothing.
\end{equation*}
Since $\rho>0$ and $l$ were arbitrary and $|n_k|\rightarrow\infty$ the point $(z,H_z)$ is recurrent in the terminology of \cite{G-H}, i.e. it is a limit point of both of the sets $\{\mathbf{T}_f^k(z,H_z):k\geq 1\}$ and $\{\mathbf{T}_f^k(z,H_z):k\leq -1\}$.
The right translation $R_{g^{-1}}:X\times G/H_z\longrightarrow X\times G/H_x$ is a homeomorphism and thus the right translate $(z,gH_x)=(z,H_z g)=R_{g^{-1}}(z,H_z)$ is also recurrent under $\mathbf{T}_f$.
We observe then that
\begin{multline*}
f(-n_k,z)\cdot g H_x=f(-n_{k'},T^{n_{k'}-n_{k}}z)\cdot f(n_{k'}-n_{k},z)\cdot g H_x \subseteq
\\ \subseteq  f(-n_{k'},T^{n_{k'}-n_{k}}z)\cdot U_{l,z} H_{z} g H_x =
f(-n_{k'},T^{n_{k'}-n_{k}}z)\cdot U_{l,z} g H_x
\end{multline*}
for all $k\geq k'$, and thus $f(-n_k,z)\cdot gH_x$ stays within a relatively compact set for all $k\geq 1$.
By passing to a subsequence we can assume that $f(-n_k,z)gH_x$ is convergent to a left coset $g'H_x$ while $T^{-n_k}z\to x$, because $T$ is an isometry.
Observe that $(x,g'H_x)$ is also in the $\mathbf{T}_f$-orbit closure of $(x,H_x)$, and from $H_x\subseteq P_x(f)$ it follows that $g'$ belongs to $P_x(f)$.
But $g'^{-1}$ belongs also to $P_x(f)$ which means that $(x,H_x)$ is in the $\mathbf{T}_f$-orbit closure of $(x,g'H_x)$, and therefore it belongs also to the orbit closure of $(z,H_z)$.
This proves that the orbit closure of $(x,H_x)$ is minimal.
From theorem 7.05 in \cite{G-H} we obtain that $(x,H_x)$ is almost periodic under $\mathbf{T}_f$ and Theorem 4.10 in \cite{G-H} then shows that this orbit closure is compact.
\end{proof}

\begin{prop}\label{pr:ER}
Let $T$ be a minimal rotation on a locally connected compact group $X$ and let $f:X\longrightarrow G$ be a continuous recurrent cocycle for $T$ with values in a nilpotent l.c.s. group $G$.
Suppose that $y\mapsto H_y$ is a consistent selection of subgroups from $X$ into $H^G$ satisfying $E_y(f)\subseteq \Stab (H_y)$ for every $y\in X$.
If there exists a point $x\in X$ with $E_x(f)=H_x$, then for every $y\in X$ there is a neighbourhood $U_y$ of $\mathbf{1}_G$ so that
\begin{equation*}
E_y(f)\cap U_y H_y=H_y.
\end{equation*}
\end{prop}

\begin{rem}
Note that if we put $H_x=\{\mathbf{1}_G\}$, then the proposition yields that equality of the sets $E_x(f)$ and $H_x$ at a single point implies that the identity $\mathbf{1}_G$ is an isolated point of $E_y(f)$ for in \emph{every} $y\in X$.
\end{rem}

In the proof of the proposition we shall need the following Lemma, which requires the same assumptions as Proposition \ref{pr:ER}:

\begin{lem}\label{lm:nsp}
Let $z\in X$ be arbitrary and assume that $g\in\Stab (H_z) \smallsetminus H_z$.
Then there exists a neighbourhood $U$ of $\mathbf{1}_G$ so that for all integers $n$
\begin{equation*}
f(n,z)\cdot g H_z\cdot f(n,z)^{-1}\cap U=f(n,z)\cdot g\cdot f(n,z)^{-1} H_{T^n z}\cap U=\varnothing.
\end{equation*}
\end{lem}

\begin{proof}
Suppose that there exists a sequence $\{(m_l,h_l)\}_{l\geq 1}\subseteq\mathbb{Z}\times H_z$ which fulfils that $f(m_l,z)\cdot g h_l\cdot f(m_l,z)^{-1}\to\mathbf{1}_G$ as $l\to\infty$.
By the compactness of $X$ we can assume that $T^{m_l}z\to y\in X$.
If we choose another integer sequence $\{k_i\}_{i\geq 1}$ with $T^{k_i}y \to z$, then the continuity of $f(k_i,\cdot)$ implies for every fixed positive integer $i$ that $f(m_l+k_i,z)\cdot g h_l\cdot f(m_l+k_i,z)^{-1}\to\mathbf{1}_G$ as $l\to\infty$.
Thus we can find an integer sequence $\{j_l\}_{l\geq 1}$ with $T^{j_l}z \to z$ and $f(j_l,z)\cdot g h_l\cdot f(j_l,z)^{-1}\to\mathbf{1}_G$ as $l\to\infty$.
Then according to the equation (\ref{eq:H_E}) and the continuity of the mapping $y\mapsto H_y$ we can select a sequence $\{(g_l,b_l)\}_{l\geq 1}\subseteq\Stab (H_z)\times G$ with $b_l\to\mathbf{1}_G$ and $f(j_l,z)=b_l  g_l$.
Hence $g_l \cdot g h_l\cdot g_l^{-1}\to\mathbf{1}_G$ as $l\to\infty$ and this contradicts the second assertion of Lemma \ref{l:double_cosets}, applied to the nilpotent l.c.s. group $\Stab (H_z)/ H_z$.
\end{proof}

\begin{proof}[Proof of Proposition \ref{pr:ER}]
We suppose that the identity coset $H_z$ is an accumulation point of $E_z(f)$ in the quotient topology of $G/H_z$ for some $z\in X$.
From $E_x(f)=H_x$ it follows for any relatively compact open neighbourhood $U_1$ of $\mathbf{1}_G$ that there exists an open neighbourhood $\mathcal{U}$ of $x$ with
\begin{equation}\label{ERy}
  E_{y}(f)\cap (\bar{U_1} H_{y}\smallsetminus U_1 H_{y})=\varnothing
\end{equation}
for all $y\in\mathcal{U}$.
Indeed, as $E_y(f)$ is a set in the quotient space $G/H_y$ for all $y\in X$ we could otherwise find a sequence $\{(y_k,g_k)\}_{k\geq 1}\subseteq X\times G$ with $y_k\to x$ and $g_k\in E_{y_k}(f)\cap (\bar{U_1} \smallsetminus U_1 H_{y_k})$.
As $\bar{U_1}$ is compact we can assume that the sequence $\{g_k\}_{k\geq 1}$ is convergent to $g'\in\bar{U_1}$.
If this limit point were in the open set $U_1 H_{x}$, then we could find an open neighbourhood $V$ of $g'$ and an open neighbourhood $W$ of some $h\in H_x$ with $VW\subseteq U_1$.
There exists a positive integer $k_0$ with $H_{y_k}\cap W\neq\varnothing$ and $g_k\in V$ for all $k\geq k_0$, because the map $y\mapsto H_y$ is continuous and $g_k\to g'$, and this contradicts that $g_k\notin U_1 H_{y_k}$ for all integers $k$.
Hence it follows that $g'\in E_x(f)\cap(\bar{U_1}\smallsetminus U_1 H_{x})$, which contradicts that $E_x(f)=H_x$.

We choose now an integer $n_1$ with $T^{n_1} z \in\mathcal{U}$ and obtain from the equations (\ref{eq:ER_ORBIT}), (\ref{eq:H_E}), and (\ref{ERy}) that the relatively compact open neighbourhood $V_1$ of $\mathbf{1}_G$ defined by $V_1=f(n_1,z)^{-1}\cdot U_1\cdot f(n_1,z)$ of $\mathbf{1}_G$ fulfils that
\begin{equation}\label{e:ER1}
  E_z(f)\cap (\bar{V_1} H_z\smallsetminus V_1 H_z)=\varnothing.
\end{equation}
The assumption that $H_z$ is not isolated in $E_z(f)$ in the quotient topology implies that there exists a $g\in\big(E_z(f)\cap V_1\big)\setminus H_z$, and Lemma \ref{lm:nsp} shows then that there exists a relatively compact open neighbourhood $U_2$ of $\mathbf{1}_G$ so that
\begin{equation*}
f(n,z)\cdot g H_z\cdot f(n,z)^{-1}\cap U_2=\varnothing
\end{equation*}
for all integers $n$.
In the same manner as above we can find an integer $n_2$ so that the neighbourhood $V_2=f(n_2,z)^{-1}\cdot U_2\cdot f(n_2,z)$ obeys that
\begin{equation}\label{e:ER2}
  E_z(f)\cap (\bar{V_2} H_z\smallsetminus V_2 H_z)=\varnothing,
\end{equation}
while on the other hand $g H_z\cap V_2=\varnothing$ and hence $g \in E_z(f)\cap (V_1 \smallsetminus \bar{V_2} H_z)$.

According to the equations (\ref{e:ER1}) and (\ref{e:ER2}) we can apply Lemma \ref{lem:eo} to find $\varepsilon>0$ so that equation (\ref{eq:eo}) holds for the neighbourhoods $V_1$ and $V_2$, and we choose then an open connected neighbourhood $\mathcal{O}$ of $z$ with $\mathcal{O}\subseteq\{y\in X:\delta(y,z)<\varepsilon/2\}$.
As the set $V_1 \smallsetminus\bar{V_2} H_z$ is a neighbourhood of $g\in E_z(f)$ we can select a sequence $\{(n_k,z_k)\}_{k\geq 1}$ in $\mathbb{Z}\times\mathcal{O}$ with $z_k\to z$, $\delta(T^{n_k} y,y)<2^{-k}\varepsilon$ for all $y\in X$, and $f(z_k,n_k)\in V_1 \smallsetminus\bar{V_2} H_{z}$.
From the connectedness of $\mathcal{O}$ and the fact that $T^{n_k}\mathcal{O}\subseteq\{y\in X:\delta(y,z)<\varepsilon\}$ it follows from Lemma \ref{lm:cont} that
\begin{equation*}
f(n_k,y)\cdot H_y\cap V_1\neq\varnothing\enspace\textup{and}\enspace f(n_k,y)\cdot H_y\cap\bar{V_2}=\varnothing
\end{equation*}
for every $y\in \mathcal{O}$, because $f(n_k,z_k)\in V_1 \smallsetminus\bar{V_2} H_z$ and both the sets $\{y\in X: f(n,y) \cdot H_y\cap V_i\neq\varnothing\}$ and $\{y\in X: f(n,y) \cdot H_y\cap \bar{V_i}=\varnothing\}$ are open.
The compactness of the set $\bar{V_1}$, the continuity of $f(n_k,\cdot)$, and the inclusion $H_y\subseteq E_y(f)$ imply in the limit $k\to\infty$ that
\begin{equation}\label{e:ERP}
  E_y(f)\cap(\bar{V_1} \smallsetminus V_2 H_y)\neq\varnothing
\end{equation}
for every $y\in\mathcal{O}$.
But as the $T$-orbit of $x$ is dense in $X$ the equation (\ref{eq:ER_ORBIT}) shows that there are points $y\in\mathcal{O}$ with $E_y(f)=H_y$, and a contradiction occurs.
\end{proof}

With these prerequisites we are now able to prove the main result on regularity of cocycles:

\begin{theo}\label{th:reg}
If $T$ is a minimal rotation on a locally connected compact group $X$ and $f:X\longrightarrow G$ is a continuous recurrent cocycle with values in a nilpotent l.c.s. group $G$, then the cocycle $f$ is regular.
\end{theo}

\begin{proof}
The proof is by induction, starting with the $0$-step nilpotent group $G=\{\mathbf{1}_G\}$ where the assertion is trivial.
We suppose now that the assertion of the theorem is fulfilled for every $(n-1)$-step nilpotent group and let $G$ be a $n$-step nilpotent group.
We denote the projection from $G$ onto $G/Z(G)$ by $\pi$ and let $\tilde{f}=\pi\circ f$ be the projection of the cocycle $f$ on the $(n-1)$-step nilpotent quotient group $G/Z(G)$, which is a regular cocycle by the induction hypothesis.

Since each of the sets $\mathcal{D}(f)$ and $\mathcal{D}\big(\tilde{f}\big)$ contains a dense $G_\delta$-set we can fix a point $x\in X$ with $E_x(f)=P_x(f)=H$ and $E_x\big(\tilde{f}\big)=P_x\big(\tilde{f}\big)=\tilde{H}$.
The Theorem \ref{t:transitive_points} shows then that $\tilde{C}=\bar{\textup{O}}_{\tilde{f}}(x,\mathbf{1}_{\tilde{G}})$ is a surjective orbit closure in $X\times\tilde{G}$ with $\tilde{H}=\Stab(\tilde{C})=\tilde{C}_x=E_x\big(\tilde{f}\big)=P_x\big(\tilde{f}\big)$ and that there is a continuous map $\tilde{\gamma}:X\longrightarrow \tilde{G}/\tilde{H}$ with $\tilde{\gamma}(y)=\tilde{C}_y=\tilde{g}_y \tilde{H}$.
Suppose that $g\in H$ and $h\in\pi^{-1}(\tilde{H})$, then there is an integer sequence $\{k_i\}_{i\geq 1}$ so that $T^{k_i}\to\textup{Id}_X$ and $f(k_i,x)\to g$ and a sequence $\{(m_l,z_l)\}_{l\geq 1}\subseteq\mathbb{Z}\times Z(G)$ with $T^{m_l}\to\textup{Id}_X$ and $f(m_l,x)\cdot z_l\to h$.
For any fixed integer $l\geq 1$ it follows that
\begin{equation*}
f(m_l,T^{k_i}x)\cdot f(k_i,x)\cdot f(-m_l,T^{m_l}x)\to f(m_l,x)\cdot z_l g z_l^{-1}\cdot f(-m_l,T^{m_l}x)
\end{equation*}
as $i\to\infty$.
In the limit $l\to\infty$ we can obtain that $hgh^{-1}\in H$ and thus $\pi^{-1}(\tilde{H})\subseteq \Stab(H)$.
Therefore the map $\gamma:X\longrightarrow G/Stab(H)$ with $y\mapsto\big(\pi^{-1}\circ\tilde{\gamma}(y)\big)\cdot\Stab(H)$ is well defined with $\gamma(x)=\Stab(H)$, and as $\pi$ is open $\gamma$ is a continuous map.
For every $y\in X$ we set $H_y=g_y H g_y^{-1}$ with some $g_y\in\gamma(y)$ and obtain a continuous map from $X$ to $H^G$ with $H_x=H$.
From $\tilde{\gamma}(T y)=\pi\big(f(y)\big)\cdot\tilde{\gamma}(y)$ for any $y\in X$ it follows that
\begin{equation*}
\gamma(Ty)=\big(\pi^{-1}\circ\tilde{\gamma}(Ty)\big)\cdot\Stab(H)=f(y)\cdot\gamma(y),
\end{equation*}
and hence the map $y\mapsto H_y$ obeys equation (\ref{eq:H_ORBIT}).
Furthermore, from $H_x=H\subseteq E_x(f)$ it follows with equations (\ref{eq:ER_ORBIT}) and (\ref{eq:H_ORBIT}) that $H_{T^kx}\subseteq E_{T^k x}(f)$ for all integers $k$.
If $y\in X$ and $g\in H_y$ are arbitrarily chosen and if $\{k_l\}_{l\geq 1}$ is a sequence with $T^{k_l}x\to y$, then by the continuity of the map $\gamma$ we can choose a sequence $\{g_l\}_{l\geq 1}\subseteq G$ with $g_l\in H_{T^{k_l}x}$ and $g_l\to g$.
We obtain that $g\in E_y(f)$ and hence both inclusions (\ref{eq:H_E}) are fulfilled.

Let $y\in X$ be arbitrary and let $g\in E_y(f)$, then it follows that $g\in\pi^{-1}\big(E_y\big(\tilde{f}\big)\big)$ and from Theorem \ref{t:essential_nil} and the inclusion $\pi^{-1}(\tilde{H})\subseteq \Stab(H)$ we can conclude that
\begin{align*}
g\in\pi^{-1}\big(\tilde{\gamma}(y)\cdot\tilde{H}\cdot\tilde{\gamma}(y)^{-1}\big)=\pi^{-1}\big(\tilde{\gamma}(y)\big)\cdot\pi^{-1}(\tilde{H})\cdot\pi^{-1}\big(\tilde{\gamma}(y)\big)^{-1}\subseteq\\
\subseteq\gamma(y)\cdot\Stab(H)\cdot\gamma(y)^{-1}=\Stab(H_y).
\end{align*}

Proposition \ref{pr:ER} shows now that for every $y\in X$ there is a neighbourhood $U_y$ of $\mathbf{1}_G$  with $E_y(f)\cap U_y H_y=H_y$ and Proposition \ref{pr:REG} applies.
Since $E_x(f)=P_x(f)=H_x$ Corollary \ref{cor:comp} shows that the cocycle $f$ is regular.
\end{proof}

In the paper \cite{A} it has been proved that the essential range of cocycle for a minimal rotation on a torus is connected if it takes its values in $\mathbb{R}^d$.
Later in the paper \cite{M1} this result has been generalised to cocycles for a minimal rotation on compact monothetic group (not necessarily connected) with values in an abelian l.c.s. group without compact subgroups.
In the next theorem we want to explore this problem in the case of nilpotent l.c.s. groups.

\begin{theo}
Let $T$ be a minimal rotation on a locally connected compact group $X$ and let $f:X\longrightarrow G$ be a continuous recurrent cocycle with values in a nilpotent l.c.s. group $G$.
Then for every $y\in X$ the essential range $E_y(f)$ is almost connected, i.e. the quotient group $E_y(f)/(E_y(f))^0$ of the essential range modulo its identity component is a compact group.
Furthermore, if $G$ is a connected and simply connected nilpotent Lie group then $E_y(f)$ is connected for every $y\in X$.
\end{theo}

\begin{proof}
Let $x$ be a point out of the set $\mathcal{D}(f)$ and let $H^0=(E_x(f))^0$ be the identity component of the essential range at $x$.
We already know from Theorem \ref{th:reg} that the cocycle $f$ is regular, and from Theorem \ref{t:essential_nil} it follows then that $E_y(f)$ is conjugate to $E_x(f)$ for all $y\in X$.
We set $H_y^0=(E_y(f))^0$ and obtain a consistent selection of subgroups from $X$ to $(H^0)^G$, because the essential ranges are conjugate on any $T$-orbit according to the equation (\ref{eq:ER_ORBIT}), while the identity component is again mapped onto identity component under conjugation.
Furthermore, the continuity of $y\mapsto H_y^0$ follows from the continuity of $y\mapsto E_y(f)$ in Theorem \ref{t:essential_nil} and the inclusion $\Stab(H)\subseteq \Stab(H^0)$.
It is well known that for any l.c.s. group $H$ the identity component $H^0$ is normal in $H$ with a totally disconnected quotient group $H/H^0$, and thus the consistent selection of subgroups $y\mapsto H_y^0$ fulfils the requirements of Proposition \ref{pr:REG}.
Now we can conclude that the $\mathbf{T}_f$-orbit closure of $(x,H^0)$ in $X\times G/H^0$ is compact, and the inclusion $H^0\subseteq P_x(f)$ implies that $\pi_{G/H^0}\big(\bar{\textup{O}}_f(x,\mathbf{1}_G)\big)=\bar{\textup{O}}_f(x,H^0)$.
Hence $E_x(f)/H^0=P_x(f)/H^0$ is compact as well as $E_y(f)/H_y^0$ for any $y\in X$, because all essential ranges are conjugate.

If $G$ is a connected and simply connected nilpotent Lie group, then $H^0$ is a connected and simply connected nilpotent Lie subgroup of $G$.
As the quotient group of any Lie group by its identity component is a discrete group, it follows that the compact quotient group $E_x(f)/H^0$ is a finite group.
If $gH^0$ is an arbitrary element in $E_x(f)/H^0$, then we can conclude that $g^k\in H^0$ for some integer $k$.
The exponential map is a diffeomorphism from the Lie Algebra $\mathfrak{g}$ onto $G$ and its restriction to the Lie algebra $\mathfrak{h}^0$ is a diffeomorphism  onto $H^0$.
If $Y$ is the unique element of $\mathfrak{g}$ with $\exp(Y)=g$, then $\exp (k Y)=\exp (Y)^k$ implies that $k Y\in\mathfrak{h}^0$, and as $\mathfrak{h}^0$ is a linear space we can conclude that $Y\in\mathfrak{h}^0$ and $g\in H^0$.
\end{proof}


\begin{thebibliography}{99}
\bibitem[A]{A} Atkinson, Giles \textit{A class of transitive cylinder transformations}, J. London Math. Soc. (2) \textbf{17} (1978), no. 2, 263--270.

\bibitem[GH]{G-H} Gottschalk, Walter Helbig; Hedlund, Gustav Arnold \textit{Topological dynamics}, American Mathematical Society Colloquium Publications, Vol. \textbf{36}, Providence, R. I., 1955.

\bibitem[LM]{LM} Lema\'nczyk, Mariusz; Mentzen, Mieczyslaw \textit{Topological ergodicity of real cocycles over minimal rotations}, Monatsh. Math. \textbf{134} (2002), no. 3, 227--246.

\bibitem[M1]{M1} Mentzen, Mieczyslaw \textit{On groups of essential values of topological cyclinder cocycles over minimal rotations}, Colloq. Math. \textbf{95} (2003), no. 2, 241--253.

\bibitem[M2]{M2} Mentzen, Mieczyslaw \textit{Some applications of groups of essential values of cocycles in topological dynamics}, Preprint, June 2004.
\end{thebibliography}
\end{document}